\documentclass[runningheads,a4paper, 10pt]{llncs}
\usepackage{amssymb, amsmath}
\usepackage{graphicx}
\usepackage{url}
\usepackage{multirow}
\newcommand{\keywords}[1]{\par\addvspace\baselineskip
\noindent\keywordname\enspace\ignorespaces#1}
\usepackage{booktabs}
\usepackage[utf8]{inputenc}

\usepackage[scale=0.9]{sourcecodepro}
\usepackage[frozencache]{minted}

\newenvironment{jlcon}
 {\VerbatimEnvironment
  \begin{minted}[
    breaklines,
    escapeinside=||,
    mathescape=true,
    frame=lines, framesep=2mm,
    fontsize=\scriptsize,
    xleftmargin=10pt, xrightmargin=10pt
    ]{jlcon}}
 {\end{minted}}

\newcommand{\julialine}[1]{\mint[fontsize=\small, breaklines, xleftmargin=10pt, xrightmargin=10pt]{julia}`#1`}
\newcommand{\juliainline}[1]{\mintinline[]{julia}`#1`}

\begin{document}

\mainmatter
\newcommand{\Z}{\mathbb{Z}}
\newcommand{\T}{\mathcal{T}}
\newcommand{\C}{\mathbb{C}}
\newcommand{\R}{\mathbb{R}}
\newcommand{\bx}{\mathbf{x}}
\newcommand{\ra}[1]{\renewcommand{\arraystretch}{#1}}
\newcommand{\polymake}{\mintinline{Perl}{polymake}}
\newcommand{\polymakejl}{\juliainline{Polymake.jl}}
\newcommand{\cxxwrapjl}{\juliainline{CxxWrap.jl}}
\newcommand{{\Cpp}}{C\nolinebreak\hspace{-.05em}\raisebox{.4ex}{\tiny\bf +}\nolinebreak\hspace{-.10em}\raisebox{.4ex}{\tiny\bf +}}

\newcommand{\TUB}{Technische Universität Berlin}
\newcommand{\CoM}{Chair of Discrete Mathematics/Geometry}
\newcommand{\AMU}{Adam Mickiewicz University in Poznań, Poland}

\title{Polymake.jl:~A new interface to {\polymake}}  
\titlerunning{Polymake.jl} 
\author{Marek Kaluba\thanks{The author was supported by the National Science Center, Poland
grant 2017/26/D/ST1/00103.
This research is carried out in the framework of the DFG funded Cluster of Excellence EXC 2046 MATH+: \textit{The Berlin Mathematics Research Center} within the Emerging Fields area.
}\inst{1,2}
\and Benjamin Lorenz\inst{1}
\and  Sascha Timme\thanks{The author was supported by the Deutsche Forschungsgemeinschaft (German Research Foundation) Graduiertenkolleg {\em Facets of Complexity} (GRK~2434).}\inst{1}
\authorrunning{Kaluba - Lorenz - Timme}
\institute{
\TUB, \CoM \\
\and
\AMU \\
}
}
\maketitle

\begin{abstract}
We present the Julia interface {\polymakejl} to {\polymake},
a software for research in polyhedral geometry.
We describe the technical design and how the integration into Julia makes it possible to combine {\polymake} with state-of-the-art numerical software.
\end{abstract}

\section{Introduction}

{\polymake} is an open source software system for computing with a wide range of objects
from polyhedral geometry and related areas \cite{polymake:2000}.
This includes convex polytopes and polyhedral fans as well as
matroids, finite permutation groups, ideals in polynomial rings and tropical varieties.
The used is interfacing the {\polymake} library through Perl language.

In this note we provide a brief overview of a new interface {\polymakejl},
which allows the use of {\polymake} in Julia \cite{julia}.
Julia is a high-level, dynamic programming language.
Distinctive aspects of Julia's design include a type system with parametric
polymorphism and multiple dispatch as its core programming paradigm.
The package {\polymakejl} can be installed, without any preparations,
using the build-in package manager that comes with Julia.
The source code is available at \url{https://github.com/oscar-system/Polymake.jl}.

\section{Functionality}
In {\polymake} the objects that a user encounters can be roughly divided into the following three classes
\begin{itemize}
  \item \emph{big objects} (e.g. cones, polytopes, simplicial complexes),
  \item \emph{small objects} (e.g. matrices, polynomials, tropical numbers),
  \item \emph{user functions}.
\end{itemize}
Broadly speaking, \emph{big objects} correspond to mathematical concepts with well defined semantics.
These can be queried, accumulate information
(e.g. a polytope defined by a set of points can ``learn'' its hyperplane representation),
and are constructed usually in Perl.
Big objects implement \emph{methods}, i.e. functions which operate on, and perform computations specific to the corresponding object.
\emph{Small objects} correspond to types or data structures which are implemented in {\Cpp}.
Standalone \emph{user functions} are exposed to the user via the Perl interpreter.

These entities are mapped to Julia in the following way:
\begin{itemize}
  \item big objects are exposed as opaque Perl objects which can be queried for their properties,
  \item small objects are wrapped through an intermediate {\Cpp} layer between Julia and \mintinline{c}{libpolymake} generated by {\cxxwrapjl},
  \item methods and user functions are mapped to Julia functions, in the case of methods, the parent object being the first argument.
\end{itemize}

A unique feature of {\polymakejl} is based on the affinity of Julia to C and {\Cpp} programming languages.
As Julia provides the possibility to call functions from dynamic libraries directly, one can call any function from the {\polymake} library as long as the function symbol is exported.
In {\polymake}, due to extensive use of templates in the {\Cpp} library, the precise definition of a function needs to be often explicitly instantiated.
Such instantiaton can be easily added to the {\polymakejl} {\Cpp} wrapper.
An example of such functionality is
\julialine{Polymake.solve_LP(inequalities, equalities, objective; sense=max)}
\noindent
function, which directly taps into {\polymake} framework for linear programming.
It is worth pointing that the signature of the exposed \juliainline{solve_LP} will
accept any instances of Julias \juliainline{AbstractMatrix} or \juliainline{AbstractVector} (where appropriate) in the paradigm of generic programming.

\section{Technical contribution}

The {\polymakejl} interface is based on {\cxxwrapjl}, a Julia package which aims
to provide a seamless interoperability between {\Cpp} and the Julia.
The interface is separated into two parts: a {\Cpp} wrapper library and a Julia package.
The former, \mintinline{c}{libpolymake} a dynamic library,
wraps the data structures (\emph{small types}) in a Julia-compatible way
and exposes functions from the callable {\Cpp} {\polymake} library.
It is then loaded through {\cxxwrapjl} where the Julia part of the package generates functions accessible from Julia.

Ihe installation of {\polymakejl} is performed through Julia package manager with the help of \juliainline{BinaryBuilder.jl} infrastructure.
Thanks to this infrastructure it is not needed for the user to perform any preparations except for installing julia itself.
All dependencies of {\polymakejl} (including the {\polymake} library, the Perl interpreter and supplementary libraries) are  installed in a binary form.
The complete installation of {\polymakejl} should take no longer than $5$ minutes on a modern hardware

Due to extensive use of metaprogramming, relatively little code was necessary to make most of the functionality of {\polymake} available in Julia:
as of version 0.3.1 {\polymakejl} consists of about $1200$ lines of {\Cpp} code and $1600$ lines of Julia code.
In particular, only the small objects need to be manually wrapped,
while functions, constructors for big objects and their methods are generated automatically from the information provided by {\polymake} itself.
This automatic code generation takes place during precompilation which is done only once during the installation.
Loading {\polymakejl} brings the familiar {\polymake} welcome banner.
\begin{jlcon}
julia> using Pkg; Pkg.add("Polymake")
  Updating registry at `~/.julia/registries/General`
  [ ... ]
  Building Polymake → `~/.julia/packages/Polymake/[...]/deps/build.log`
julia> using Polymake
[ Info: Precompiling Polymake [d720cf60-89b5-51f5-aff5-213f193123e7]
[ Info: Generating module common
[ Info: Generating module ideal
[ Info: Generating module graph
[ Info: Generating module fulton
[ Info: Generating module fan
[ Info: Generating module group
[ Info: Generating module polytope
[ Info: Generating module topaz
[ Info: Generating module tropical
[ Info: Generating module matroid

polymake version 4.0
Copyright (c) 1997-2020
Ewgenij Gawrilow, Michael Joswig, and the polymake team
Technische Universität Berlin, Germany
https://polymake.org

This is free software licensed under GPL; see the source for
copying conditions.
There is NO warranty; not even for MERCHANTABILITY or
FITNESS FOR A PARTICULAR PURPOSE.

\end{jlcon}
  
\subsection*{Big Objects}  

All big objects are constructed by direct calls to their constructors, e.g.
\julialine{polytope.Polytope(POINTS=[1 1 2; 1 3 4])}
\noindent
constructs a rational polytope from (homogeneous) coordinates of points given row-wise.
We attach the {\polymake} docstring to the structure such that the documentation is readily available in Julia.
\begin{jlcon}
help?> polytope.Polytope
Not necessarily bounded convex polyhedron, i.e., the feasible region of a linear program.
[...]
\end{jlcon}
Template parameters can also be passed to big objects, e.g., to construct a polytope with floating point precision it is sufficient to call
\julialine{polytope.Polytope{Float64}(...) .}
\noindent A caveat is that all parameters must be valid Julia objects.\footnote{For advanced use (when this is not the case) we provide the \texttt{@pm} macro.}
The properties of big objects are accessible through the \juliainline{bigobject.property} syntax which mirrors the \mintinline{Perl}{$bigobject->property} syntax in {\polymake}.
Note that some properties of big objects in {\polymake} are indeed methods with no arguments and therefore in Julia they are only available as such.

\subsection*{Small objects}
The list of small objects available in {\polymakejl} includes basic types such as 
arbitrary size integers (subtypes \juliainline{Integer}),
rationals (subtypes \juliainline{Real}),
vectors and matrices (subtype \juliainline{AbstractArray}), and many more.
These data types can be converted to appropriate Julia types, but are also subtypes of the corresponding Julia abstract types (as indicated above).
This allows to use {\polymakejl} types in generic methods, which is the paradigm of julia programming.

As already mentioned, these small objects need to be manually wrapped in the {\Cpp} part of {\polymakejl}.
In particular, all possible combinations of such types, e.g., an array of sets of rationals, need to be explicitly wrapped.
Note that {\polymake} is able to generate dynamically any combination of small objects.
Thus, we cannot guarantee that all small objects a user will encounter is covered.
However, the small objects available in {\polymakejl} are sufficient for the most common use cases.

\subsection*{Functions}

A function in {\polymakejl} calling {\polymake} may return either a big or a small object,
and the generic return type (\juliainline{PropertyValue}, a container opaque to Julia) is transparently converted to one of the known (small) data types.
\footnote{This conversion can be deactivated by adding \juliainline{PropertyValue} as the first argument to function/method call.}
If the data type of the returned function value is not known to {\polymakejl}, the conversion fails and an instance of \juliainline{PropertyValue} is returned.
It can be either passed back as an argument to a {\polymakejl} function, or converted to a known type using the \juliainline{@convert_to} macro.

\begin{jlcon}
julia> K5 = graph.complete(5);
julia> K5.MAX_CLIQUES
PropertyValue wrapping pm::PowerSet<long, pm::operations::cmp>
{{0 1 2 3 4}}
julia> @convert_to Array{Set} K5.MAX_CLIQUES
pm::Array<pm::Set<long, pm::operations::cmp> >
{0 1 2 3 4} 
\end{jlcon}

All user functions from {\polymake} are available in modules corresponding to their applications, e.g.
\juliainline{homology} functions from the application \juliainline{topaz} can be called as \juliainline{topaz.homology(...)} in Julia.
Moreover {\polymake} docstrings for functions are available in Julia to allow for easy help\footnote{The documentation currently uses the Perl syntax}:
\begin{jlcon}
julia> ?topaz.homology
  homology(complex, co; Options)
  
   Calculate the reduced (co-)homology groups of a simplicial complex.
  
  Arguments:
    Array<Set<Int>> complex 
    Bool co set to true for cohomology
  
  Options: 
    dim_low => Int narrows the dimension range of interest, with negative values being treated as co-dimensions
    dim_high => Int see dim_low

    [ ... ]
\end{jlcon}

\paragraph*{Function Arguments}

Functions in {\polymakejl} accept as their arguments: simple data types (bools, machine integers, floats), wrapped native types, or objects returned by {\polymake} (e.g. \juliainline{BigObject}, or \juliainline{PropertyValue}).
Due to the easy extendability of methods in Julia, a foreign type could be passed seamlessly to {\polymakejl} function if an appropriate \juliainline{Base.convert} method, which return one of the above types, is defined:
\julialine{Base.convert(::Type{Polymake.PolymakeType}, x::ForeignType)}

\texttt{Polymake.jl} also wraps the extensive visualization methods of {\polymake} which can be used to produce images and animations of geometric objects.
These include the interactive visualizations using \texttt{three.js}.
Due to the convenient extendability of Julia the visualization also integrates seamlessly with Jupyter notebooks.

\section{Example}

This section demonstrates the interface of {\polymakejl} on a concrete example.
An advantage the package is allow effortless combination of computations in polyhedral geometry with e.g. state-of-the-art numerical software.
Here we combine{\polymakejl} with \juliainline{HomotopyContinuation.jl} \cite{HomotopyContinuation.jl}, a Julia package for the numerical solution of systems of polynomial equations.
In particular, we test a theoretical result from Soprunova and Sottile \cite{Soprunova:Sottile:2006} on non-trivial lower bounds for the number of real solutions to sparse polynomial systems.

The results show how we can construct a sparse polynomial system that has a non-trivial lower bound on the number of real solutions starting from an integral point configuration.
We start with the 10 lattice points $A=\{a_1,\ldots,a_{10}\} \subset \Z^2$ of the scaled two-dimensional simplex $3\Delta_2$ and look at the regular triangulation $\T$ induced by the lifting $\lambda = [12, 3, 0, 0, 8, 1, 0, 9, 5, 15]$.
\begin{jlcon}
julia> A = polytope.lattice_points(polytope.simplex(2,3));
julia> |$\lambda$| = [12, 3, 0, 0, 8, 1, 0, 9, 5, 15];
julia> F = polytope.regular_subdivision(A, |$\lambda$|);
julia> T = topaz.GeometricSimplicialComplex(COORDINATES = A[:,2:end], FACETS = F)
type: GeometricSimplicialComplex<Rational>

COORDINATES
        0 0
        0 1
        [ ... ]

FACETS
  pm::Set<long, pm::operations::cmp>
{5 6 8}
  pm::Set<long, pm::operations::cmp>
{5 7 8}
  [ ... ]
\end{jlcon}

\begin{figure}[h]
  \begin{center}
      \includegraphics[width=0.3\textwidth]{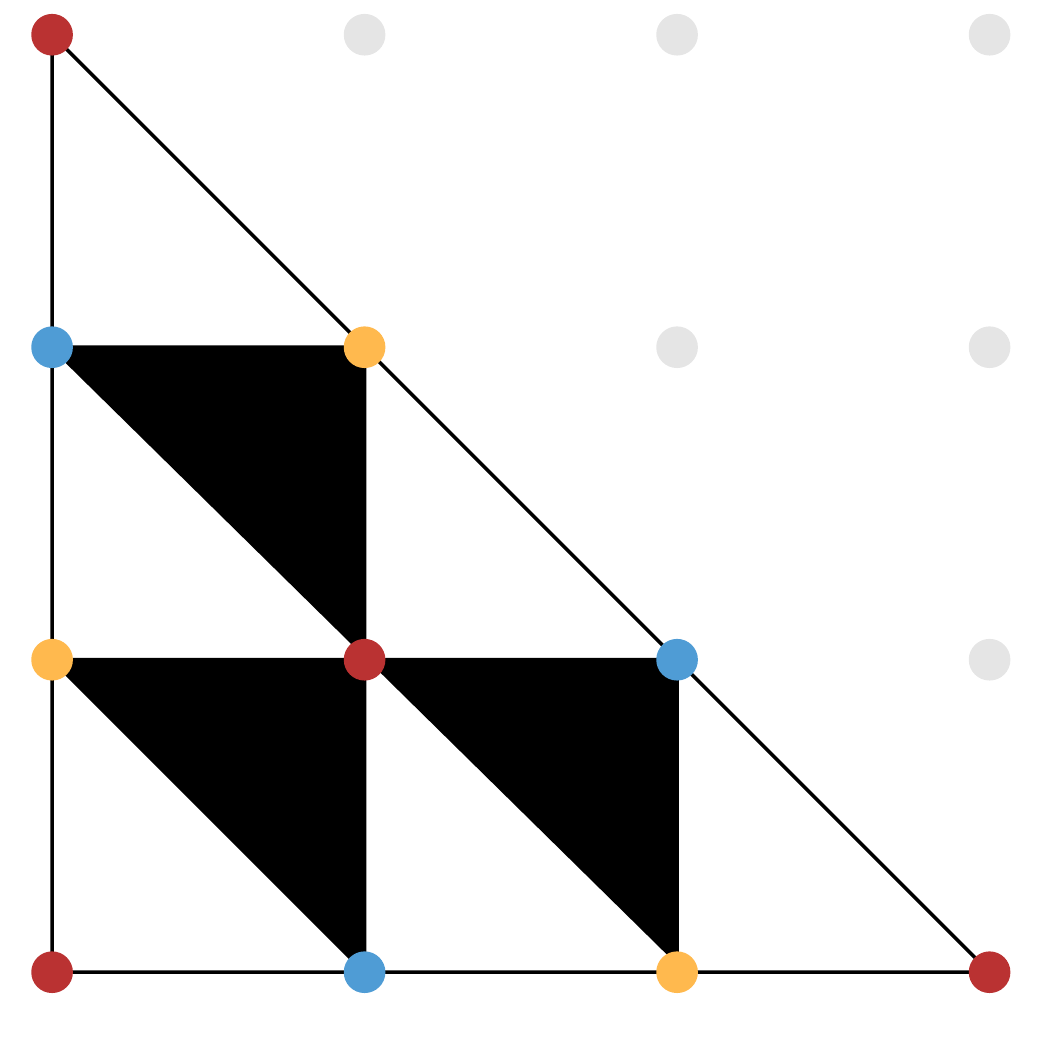}
  \end{center}
  \caption{Foldable subdivision of $3\Delta_2$.}\label{fig:triangulation}
\end{figure}

\noindent The triangulation $\T$ is very special in that it is \emph{foldable} (or ``balanced''), i.e., the dual graph is bipartite. 
This means that the triangles can be colored, say, black and white such that no two triangles of the same color share an edge.
See Figure \ref{fig:triangulation} for an illustration.
The \emph{signature} $\sigma(\T)$ of a balanced triangulation of a polygon is the absolute value of the difference of the number of black triangles and the number of the white triangles whose normalized volume is odd.
The vertices of a foldable triangulation can be colored by $d+1$ colors \cite{Joswig:2002} (such that vertices of the same color do not share an edge), where $d$ is the dimension. Here $d=2$, so $3$ colors suffice.
We can check both properties with {\polymake}.
\begin{jlcon}
julia> (foldable = T.FOLDABLE, signature = T.SIGNATURE)
  (foldable = true, signature = 3)
\end{jlcon}

Now, a \emph{Wroński polynomial} $W_{\T,s}(x)$ has the lifted lattice points as exponents, and only one non-zero coefficient $c_i \in \R$ per color class of vertices of the triangulation
$$W_{\T,s}(x) = \sum_{i=0}^d c_i \left (\sum_{j:\;\textnormal{\scriptsize color($a_j$)=$i$}} s^{\lambda_i} x^{a_j} \right ) \,.$$
A \emph{Wroński system} consists of $d$ Wroński polynomials with respect to the same lattice points $A$ and lifting $\lambda$
such that for general $s = s_0 \in [0,1]$ it has precisely $d!\operatorname{vol}(\operatorname{conv}(A))$ distinct complex solutions,
which is the highest possible number by Kushnirenko’s Theorem \cite{Kushnirenko:1975}.

Soprunova and Sottile showed that a Wroński system has at least $\sigma(\T)$ distinct real solutions
if two conditions are satisfied.
First, a certain double cover of the real toric variety associated with $A$ must be orientable. This is the case here.
Second, the \emph{Wroński center ideal}, a zero-dimensional ideal in coordinates $x_1,x_2$ and $s$ depending on $\T$, has no real roots with $s$ coordinate between 0 and 1.
Let us verify this condition using \juliainline{HomotopyContinuation.jl}.
{\polymake} already has an implementation of the Wroński center ideal for us to use.
However, we have to convert the ideal returned by {\polymakejl} to a polynomial system which \juliainline{HomotopyContinuation.jl} understands.
This can be accomplished with a simple routine.
\begin{jlcon}
julia> using HomotopyContinuation
julia> function hc_poly(f, vars)
          M = Polymake.monomials_as_matrix(f)
          monomials = [prod(vars.^m) for m in eachrow(M)]
          coeffs = Int.(Polymake.coefficients_as_vector(f))
          sum(map(*, coeffs, monomials))
      end;
julia> I = polytope.wronski_center_ideal(A, |$\lambda$|)
julia> @polyvar x[1:2] s;
julia> HC_I = [hc_poly(f, [x;s]) for f in I.GENERATORS]
  3-element Array{Polynomial{true,Int},1}:
    |$x_1^3s^{15} + s^{12} + x_1 x_2 s + x_2^3$|
    |$x_1^2s^{9} + x_2s^{3} + x_1 x_2^2$|
    |$x_1s^{8} + x_1^2 x_2 s^{5} + x_2^2$|
\end{jlcon}
Since we are only interested in solutions in the algebraic torus $(\mathbb{C}^*)^3$
we can use polyhedral homotopy \cite{Huber:Sturmfels:1995} to efficiently compute the solutions.
\begin{jlcon}
julia> @time res = solve(HC_I; start_system = :polyhedral, only_torus = true)
  0.010595 seconds (3.03 k allocations: 215.766 KiB)
Result{Array{Complex{Float64},1}} with 54 solutions
===================================================
* 54 non-singular solutions (2 real)
* 0 singular solutions (0 real)
* 54 paths tracked
* random seed: 782949
\end{jlcon}
Out of the 54 complex roots only two solutions are real.
By closer inspection, we see that no solution has the $s$-coordinate in $(0,1)$.
\begin{jlcon}
julia> HomotopyContinuation.real_solutions(res)
2-element Array{Array{Float64,1},1}:
 [-0.2117580095433453, -215.72260079314424, 4.411470567441922]   
 [-0.6943590430596768, -0.41424188458258815, -0.8952189506082179]
\end{jlcon}
\noindent Therefore, the Wroński system with respect to $A$ and $\lambda$ for $s=1$ has at least $\sigma(\T)=3$ real solutions. Let us verify this on an example.
\begin{jlcon}
julia> c = Vector{Polymake.Rational}[[19,8,-19], [39,7,42]];
julia> W = polytope.wronski_system(A, |$\lambda$|, c, 1)
julia> HC_W = [hc_poly(f, x) for f in W.GENERATORS];
julia> W_res = HomotopyContinuation.solve(HC_W)
Result{Array{Complex{Float64},1}} with 9 solutions
==================================================
* 9 non-singular solutions (3 real)
* 0 singular solutions (0 real)
* 9 paths tracked
* random seed: 813729
\end{jlcon}
\noindent Finally, we can use the \juliainline{ImplicitPlots.jl} package to visualize the real solutions of the Wroński system \juliainline{W}.
\begin{jlcon}
  julia> W_real = HomotopyContinuation.real_solutions(W_res)
  julia> using ImplicitPlots, Plots;
  julia> p = plot(aspect_ratio = :equal);
  julia> implicit_plot!(p, HC_W[1]);
  julia> implicit_plot!(p, HC_W[2]; color=:indianred);
  julia> scatter!(first.(r), last.(r), markercolor=:black)
\end{jlcon}
\begin{figure}[h]
  \begin{center}
      \includegraphics[width=0.4\textwidth]{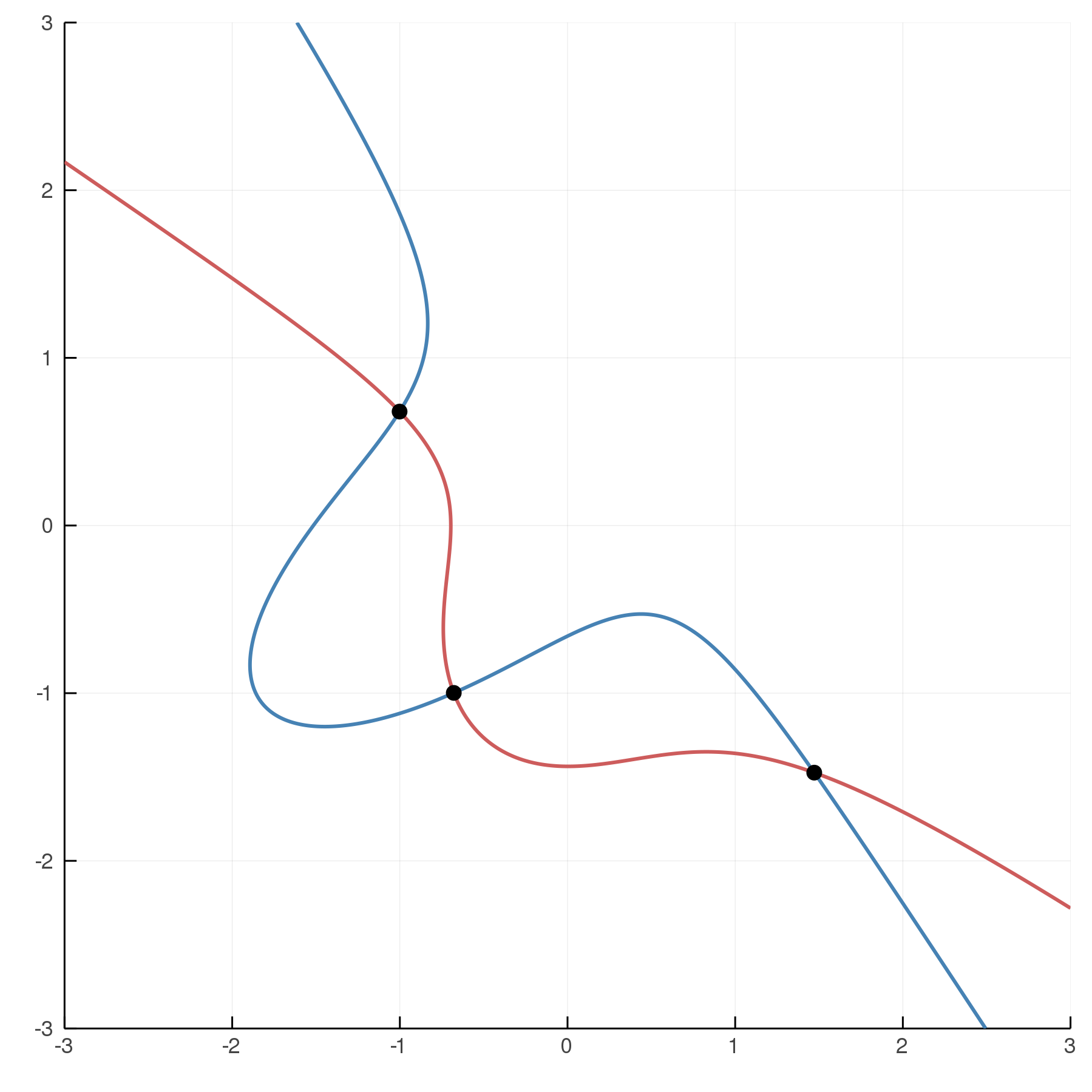}
    \end{center}
    \caption{Visualization of the Wroński system $W$ and its $3$ solutions.}
\end{figure}

\section*{Acknowledgements}
We would like to express our thanks to Alexej Jordan and Sebastian Gutsche for all their help during the development of {\polymakejl}.


\end{document}